\documentclass[a4paper,12pt,oneside]{article}
\usepackage{amsmath,amsfonts,amssymb,amsmath,latexsym,amsthm,enumerate}
\usepackage{hyperref}

\newtheorem{thm}{Theorem}

\theoremstyle{definition}

\theoremstyle{plain}

\begin{document}
\title {Poly-Cauchy numbers and polynomials of the second kind}
\author{by \\Dae San Kim and Taekyun Kim}\date{}\maketitle

\begin{abstract}
\noindent In this paper, we consider the poly-cauchy polynomials and numbers of the second kind which were studied by Komatsu in $\lbrack 10\rbrack$. We note that the poly-Cauchy polynomials of the second kind are the special generalized Bernoulli polynomials of the second kind. The purpose of this paper is to give various identities of the poly-Cauchy polynomials of the second kind which are derived from umbral calculus.
\end{abstract}

\section{Introduction}

As is well  known, the Bernoulli polynomials of the second kind are defined by the generating function to be
\begin{equation}\label{eq:1}
\frac{t}{\log{(1+t)}}(1+t)^{x}=\sum_{n=0}^{\infty}b_{n}(x)\frac{t^{n}}{n!},\,\,\,\,(\text{see}\,\,\lbrack 12\rbrack).
\end{equation}
When $x=0$, $b_{n}=b_{n}(0)$ are called the Bernoulli numbers of the second kind.\\
Let $Lif_{k}(x)$ be the polylogarithm factorial function which is defined by
\begin{equation}\label{eq:2}
Lif_{k}(x)=\sum_{n=0}^{\infty}\frac{x^{m}}{m!(m+1)^{k}},\,\,\,\,(\text{see}\,\,\lbrack 9-11\rbrack).
\end{equation}
The poly-Cauchy polynomials of the second kind $\tilde{C}_{n}^{(k)}(x)$ $(k\in\mathbf{Z},\,\,n\in\mathbf{Z}_{\geq 0})$ are defined by the generating function to be
\begin{equation}\label{eq:3}
Lif_{k}\left(-\log{(1+t)}\right)(1+t)^{x}=\sum_{n=0}^{\infty}\tilde{C}_{n}^{(k)}(x)\frac{t^{n}}{n!},\,\,\,\,(\text{see}\,\,\lbrack 10\rbrack).
\end{equation}
When $x=0$, $\tilde{C}_{n}^{(k)}=\tilde{C}_{n}^{(k)}(0)$ are called the poly-Cauchy numbers.\\
In particular, if we take $k=1$, then we have
\begin{equation}\label{eq:4}
Lif_{1}\left(-\log{(1+t)}\right)(1+t)^{x}=\frac{t}{(1+t)\log{(1+t)}}(1+t)^{x}=\frac{t(1+t)^{x-1}}{\log{(1+t)}}.
\end{equation}
Thus, we note that
\begin{equation}\label{eq:5}
\tilde{C}_{n}^{(1)}(x)=b_{n}(x-1)=B_{n}^{(n)}(x),
\end{equation}
where $B_{n}^{(\alpha)}(x)$ are called the Bernoulli polynomials of order $\alpha$ (see $\lbrack 3,4,5,6,7\rbrack$).\\
When $x=0$, $\tilde{C}_{n}^{(1)}=\tilde{C}_{n}^{(1)}(0)=b_{n}(-1)=B_{n}^{(n)}$, where $B_{n}^{(\alpha)}$ are called the Bernoulli numbers of order $\alpha$.\\
The falling factorial is defined by
\begin{equation}\label{eq:6}
(x)_{n}=x(x-1)\cdots(x-n+1)=\sum_{l=0}^{n}S_{1}(n,l)x^{l},
\end{equation}
where $S_{1}(n,l)$ is the stirling number of the first kind.\\
For $m\in\mathbf{Z}_{\geq 0}$, it is well known that
\begin{align}\label{eq:7}
\left(\log(1+t)\right)^{m}&=m!\sum_{l=m}^{\infty}S_{1}(l,m)\frac{t^{l}}{l!}\\
&=\sum_{l=0}^{\infty}S_{1}(l+m,m)\frac{m!}{(l+m)!}t^{l+m},\,\,\,\,(\text{see}\,\,\lbrack 12,13\rbrack).\nonumber
\end{align}
For $\lambda\in\mathbf{C}$ with $\lambda\neq 1$, the Frobenius-Euler polynomials of order $r$ are defined by the generating function to be
\begin{equation*}
\left(\frac{1-\lambda}{e^{t}-\lambda}\right)^{r}e^{xt}=\sum_{n=0}^{\infty}H_{n}^{(r)}\left(x\vert\lambda\right)\frac{t^{n}}{n!},\,\,\,\,(\text{see}\,\,\lbrack 1,3,7,8\rbrack).
\end{equation*}
Let $\mathbf{C}$ be the complex number field and let $\mathcal{F}$ be the set of all formal power series in the variable $t$:
\begin{equation}\label{eq:8}
\mathcal{F}=\left\{f(t)=\sum_{k=0}^{\infty}\frac{a_{k}}{k!}t^{k}\Bigg\vert a_{k}\in\mathbf{C}\right\}.
\end{equation}
Let $\mathbb{P}=\mathbf{C}\lbrack x\rbrack$ and let $\mathbb{P}^{*}$ be the vector space of all linear functionals on $\mathbb{P}$. $\langle L\vert p(x)\rangle$ is the action of the linear functional $L$ on the polynomial $p(x)$, and we recall that the vector space operations on $\mathbb{P}^{*}$ are defined by $\langle L+M\vert p(x)\rangle=\langle L\vert p(x)\rangle +\langle M\vert p(x)\rangle$, $\langle cL\vert p(x)\rangle=c\langle L\vert p(x)\rangle$, where $c$ is a complex constant in $\mathbf{C}$. For $f(t)\in\mathcal{F}$, let us define the linear functional on $\mathbb{P}$ by setting
\begin{equation}\label{eq:9}
\langle f(t)\vert x^{n}\rangle=a_{n},\,\,\,\,(n\geq 0).
\end{equation}
Then, by (\ref{eq:8}) and (\ref{eq:9}), we get
\begin{equation}\label{eq:10}
\langle t^{k}\vert x^{n}\rangle=n!\delta_{n,k},\,\,\,\,(n,k\geq 0),
\end{equation}
where $\delta_{n,k}$ is the Kronecker's symbol.\\
For $f_{L}(t)=\sum_{k=0}^{\infty}\frac{\langle L\vert x^{k}\rangle}{k!}t^{k}$, we have $\langle f_{L}(t)\vert x^{n}\rangle=\langle L\vert x^{n}\rangle$. That is, $L=f_{L}(t)$. The map $L\longmapsto f_{L}(t)$ is a vector space isomorphism from $\mathbb{P}^{*}$ onto $\mathcal{F}$. Henceforth, $\mathcal{F}$ denotes both the algebra of formal power series in $t$ and the vector space of all linear functionals on $\mathbb{P}$, and so an element $f(t)$ of $\mathcal{F}$ will be thought of as both a formal power series and a linear functional. We call $\mathcal{F}$ the umbral algebra and the umbral calculus is the study of umbral algebra. The order $O\left(f(t)\right)$ of a power series $f(t)(\neq 0)$ is the smallest integer $k$ for which the coefficient of $t^{k}$ does not vanish. If $O\left(f(t)\right)=1$, then $f(t)$ is called a delta series ; if $O\left(f(t)\right)=0$, then $f(t)$ is called an invertible series (see $\lbrack 2,12,13\rbrack$). For $f(t), g(t)\in\mathcal{F}$ with $O\left(f(t)\right)=1$ and $O\left(g(t)\right)=0$, there  exists a unique sequence $s_{n}(x)$ ($\deg{s_{n}(x)}=n$) such that $\langle g(t)f(t)^{k}\vert s_{n}(x)\rangle=n!\delta_{n,k}$ for $n, k\geq 0$. The sequence $s_{n}(x)$ is called the Sheffer sequence for $\left(g(t), f(t)\right)$ which is denoted by $s_{n}(x)\sim\left(g(t), f(t)\right)$ (see $\lbrack 12,13\rbrack$).\\
For $f(t), g(t)\in\mathcal{F}$ and $p(x)\in\mathbb{P}$, we have
\begin{equation}\label{eq:11}
\langle f(t)g(t)\vert p(x)\rangle=\langle f(t)\vert g(t)p(x)\rangle=\langle g(t)\vert f(t)p(x)\rangle,
\end{equation}
and
\begin{equation}\label{eq:12}
f(t)=\sum_{k=0}^{\infty}\langle f(t)\vert x^{k}\rangle\frac{t^{k}}{k!},\,\,\,\, p(x)=\sum_{k=0}^{\infty}\langle t^{k}\vert p(x)\rangle\frac{x^{k}}{k!}.
\end{equation}
Thus, by (\ref{eq:12}), we get
\begin{equation}\label{eq:13}
t^{k}p(x)=p^{(k)}(x)=\frac{d^{k}p(x)}{dx^{k}},\,\,\,\,\text{and}\,\,\,\, e^{yt}p(x)=p(x+y).
\end{equation}
Let us assume that $s_{n}(x)\sim\left(g(t),f(t)\right)$. Then the generating function of $s_{n}(x)$ is given by
\begin{equation}\label{eq:14}
\frac{1}{g\left(\bar{f}(t)\right)}e^{x\bar{f}(t)}=\sum_{n=0}^{\infty}s_{n}(x)\frac{t^{n}}{n!},\,\,\,\,\text{for all}\,\, x\in\mathbf{C},
\end{equation}
where $\bar{f}(t)$ is the compositional inverse of $f(t)$ with $\bar{f}\left(f(t)\right)=t$ (see $\lbrack 12,13\rbrack$).\\
For $s_{n}(x)\sim\left(g(t),f(t)\right)$, we have the following equation:
\begin{equation}\label{eq:15}
f(t)s_{n}(x)=ns_{n-1}(x),\,\,\,\,(n\geq 0),
\end{equation}
\begin{equation}\label{eq:16}
s_{n}(x)=\sum_{j=0}^{n}\frac{1}{j!}\left\langle g\left(\bar{f}(t)\right)^{-1}\bar{f}(t)^{j}\big\vert x^{n}\right\rangle x^{j},
\end{equation}
and
\begin{equation}\label{eq:17}
s_{n}(x+y)=\sum_{j=0}^{n}\binom{n}{j}s_{j}(x)p_{n-j}(y),
\end{equation}
where $p_{n}(x)=g(t)s_{n}(x)$, (see $\lbrack 12,13\rbrack$).\\
Let us assume that $p_{n}(x)\sim\left(1,f(t)\right)$, $q_{n}(x)\sim\left(1,g(t)\right)$. Then the transfer formula is given by
\begin{equation*}
q_{n}(x)=x\left(\frac{f(t)}{g(t)}\right)^{n}x^{-1}p_{n}(x),\,\,\,\,(n\geq 0),\,\,(\text{see}\,\,\lbrack 12\rbrack).
\end{equation*}
For $s_{n}(x)\sim\left(g(t),f(t)\right)$, $r_{n}(x)\sim\left(h(t),l(t)\right)$, let us assume that
\begin{equation}\label{eq:18}
s_{n}(x)=\sum_{m=0}^{n}C_{n,m}r_{n}(x),\,\,\,\,(n\geq 0).
\end{equation}
Then we have
\begin{equation}\label{eq:19}
C_{n,m}=\frac{1}{m!}\left\langle\frac{h\left(\bar{f}(t)\right)}{g\left(\bar{f}(t)\right)}l\left(\bar{f}(t)\right)^{m}\Bigg\vert x^{n}\right\rangle,\,\,\,\,(\text{see}\,\,\lbrack 12\rbrack).
\end{equation}
In this paper, we investigate the properties of the poly-Cauchy numbers and polynomials of the second kind with umbral calculus viewpoint. The purpose of this paper is to give various identities of the poly-Cauchy polynomials of the second kind which are derived from umbral calculus.

\section{Poly-Cauchy numbers and polynomials of the second kind}

From (\ref{eq:3}), we note that $\tilde{C}_{n}^{(k)}(x)$ is the Sheffer sequence for the pair
\begin{equation*}
\left(g(t)=\frac{1}{Li{f_{k}}(-t)},f(t)=e^{t}-1\right),
\end{equation*}
that is,
\begin{equation}\label{eq:20}
\tilde{C}_{n}^{(k)}(x)\sim\left(\frac{1}{Li{f_{k}}(-t)},e^{t}-1\right).
\end{equation}
Thus, by (\ref{eq:14}) and (\ref{eq:20}), we get
\begin{equation}\label{eq:21}
\sum_{n=0}^{\infty}\tilde{C}_{n}^{(k)}(x)\frac{t^{n}}{n!}=Li{f_{k}}\left(-\log{(1+t)}\right)(1+t)^{x},
\end{equation}
and
\begin{equation}\label{eq:22}
\sum_{n=0}^{\infty}\tilde{C}_{n}^{(k)}\frac{t^{n}}{n!}=Li{f_{k}}\left(-\log{(1+t)}\right).
\end{equation}
From (\ref{eq:20}), we have
\begin{equation}\label{eq:23}
\frac{1}{Li{f_{k}}(-t)}\tilde{C}_{n}^{(k)}(x)\sim\left(1,e^{t}-1\right),
\end{equation}
and
\begin{equation}\label{eq:24}
(x)_{n}=\sum_{l=0}^{n}S_{1}(n,l)x^{l}\sim\left(1,e^{t}-1\right).
\end{equation}
By (\ref{eq:23}) and (\ref{eq:24}), we get
\begin{align}\label{eq:25}
\tilde{C}_{n}^{(k)}(x)&=Li{f_{k}}(-t)(x)_{n}=\sum_{m=0}^{n}S_{1}(n,m)Li{f_{k}}(-t)x^{m}\\
&=\sum_{m=0}^{n}S_{1}(n,m)\sum_{a=0}^{m}\frac{(-1)^{a}}{a!(a+1)^{k}}t^{a}x^{m}\nonumber\\
&=\sum_{m=0}^{n}\sum_{a=0}^{m}S_{1}(n,m)\frac{(-1)^{a}\binom{m}{a}}{(a+1)^{k}}x^{m-a}\nonumber\\
&=\sum_{m=0}^{n}\sum_{j=0}^{m}S_{1}(n,m)\frac{(-1)^{m-j}\binom{m}{j}}{(m-j+1)^{k}}x^{j}\nonumber\\
&=\sum_{j=0}^{n}\left\{\sum_{m=j}^{n}S_{1}(n,m)\frac{(-1)^{m-j}\binom{m}{j}}{(m-j+1)^{k}}\right\}x^{j}.\nonumber
\end{align}
By (\ref{eq:16}) and (\ref{eq:20}), we get
\begin{equation}\label{eq:26}
\tilde{C}_{n}^{(k)}(x)=\sum_{j=0}^{n}\frac{1}{j!}\left\langle Li{f_{k}}\left(-\log{(1+t)}\right)\left(\log{(1+t)}\right)^{j}\Big\vert x^{n}\right\rangle x^{j}.
\end{equation}
Now, we observe that
\begin{align}\label{eq:27}
&\left\langle Li{f_{k}}\left(-\log{(1+t)}\right)\left(\log{(1+t)}\right)^{j}\Big\vert x^{n}\right\rangle\\
&=\sum_{m=0}^{\infty}\frac{(-1)^{m}}{m!(m+1)^{k}}\left\langle\left(\log{(1+t)}\right)^{m+j}\Big\vert x^{n}\right\rangle\nonumber\\
&=\sum_{m=0}^{n-j}\frac{(-1)^{m}}{m!(m+1)^{k}}\sum_{l=0}^{n-j-m}\frac{S_{1}(l+m+j,m+j)}{(l+m+j)!}(m+j)!\left\langle t^{m+j+l}\big\vert x^{n}\right\rangle\nonumber\\
&=\sum_{m=0}^{n-j}\frac{(-1)^{m}}{m!(m+1)^{k}}\sum_{l=0}^{n-m-j}\frac{S_{1}(l+m+j,m+j)}{(l+m+j)!}(m+j)!n!\delta_{n,l+m+j}\nonumber\\
&=\sum_{m=0}^{n-j}\frac{(-1)^{m}(m+j)!}{m!(m+1)^{k}}S_{1}(n,m+j).\nonumber
\end{align}
From (\ref{eq:26}) and (\ref{eq:27}), we have
\begin{align}\label{eq:28}
\tilde{C}_{n}^{(k)}(x)&=\sum_{j=0}^{n}\frac{1}{j!}\sum_{m=0}^{n-j}\frac{(-1)^{m}(m+j)!}{m!(m+1)^{k}}S_{1}(n,m+j)x^{j}\\
&=\sum_{j=0}^{n}\left\{\sum_{m=0}^{n-j}\frac{(-1)^{m}\binom{m+j}{m}}{(m+1)^{k}}S_{1}(n,m+j)\right\}x^{j}\nonumber\\
&=\sum_{j=0}^{n}\left\{\sum_{m=j}^{n}\frac{(-1)^{m-j}\binom{m}{j}}{(m-j+1)^{k}}S_{1}(n,m)\right\}x^{j}.\nonumber
\end{align}
From (\ref{eq:1}), we note that
\begin{equation}\label{eq:29}
\frac{1}{Li{f_{k}}(-t)}\tilde{C}_{n}^{(k)}(x)\sim\left(1,e^{t}-1\right),\,\,\,\,x^{n}\sim(1,t).
\end{equation}
For $n\geq 1$, by (\ref{eq:18}) and (\ref{eq:29}), we get
\begin{align}\label{eq:30}
\frac{1}{Li{f_{k}}(-t)}\tilde{C}_{n}^{(k)}(x)&=x\left(\frac{t}{e^{t}-1}\right)^{n}x^{-1}x^{n}=x\left(\frac{t}{e^{t}-1}\right)^{n}x^{n-1}\\
&=xB_{n-1}^{(n)}(x)=\sum_{l=0}^{n-1}\binom{n-1}{l}B_{n-1-l}^{(n)}x^{l+1}.\nonumber
\end{align}
Thus, by (\ref{eq:30}), we see that
\begin{align}\label{eq:31}
\tilde{C}_{n}^{(k)}(x)&=\sum_{l=0}^{n-1}\binom{n-1}{l}B_{n-1-l}^{(n)}Li{f_{k}}(-t)x^{l+1}\\
&=\sum_{l=0}^{n-1}\sum_{m=0}^{l+1}(-1)^{m}\binom{n-1}{l}\binom{l+1}{m}\frac{B_{n-1-l}^{(n)}}{(m+1)^{k}}x^{l+1-m}\nonumber\\
&=\sum_{l=0}^{n-1}\sum_{j=0}^{l+1}(-1)^{l+1-j}\binom{n-1}{l}\binom{l+1}{j}\frac{B_{n-1-l}^{(n)}}{(l+2-j)^{k}}x^{j}\nonumber\\
&=\sum_{l=0}^{n-1}(-1)^{l+1}\binom{n-1}{l}\frac{B_{n-1-l}^{(n)}}{(l+2)^{k}}\nonumber\\
&\quad+\sum_{j=1}^{n}\left\{\sum_{l=j-1}^{n-1}(-1)^{l+1-j}\binom{n-1}{l}\binom{l+1}{j}\frac{B_{n-1-l}^{(n)}}{(l+2-j)^{k}}\right\}x^{j}.\nonumber
\end{align}
Therefore, by (\ref{eq:25})(or (\ref{eq:28})) and (\ref{eq:31}), we obtain the following theorem.

\begin{thm}\label{eq:thm1}
For $n\geq 1$, $1\leq j\leq n$, we have
\begin{equation*}
\sum_{m=j}^{n}\frac{(-1)^{m-j}\binom{m}{j}}{(m-j+1)^{k}}S_{1}(n,m)=\sum_{l=j-1}^{n-1}(-1)^{l+1-j}\binom{n-1}{l}\binom{l+1}{j}\frac{B_{n-1-l}^{(n)}}{(l+2-j)^{k}}.
\end{equation*}
In addition, for $n\geq 1$, we have
\begin{align*}
\tilde{C}_{n}^{(k)}&=\sum_{m=0}^{n}S_{1}(n,m)\frac{(-1)^{m}}{(m+1)^{k}}\\
&=\sum_{l=0}^{n-1}(-1)^{l+1}\binom{n-1}{l}\frac{B_{n-1-l}^{(n)}}{(l+2)^{k}}.
\end{align*}
\end{thm}

\noindent From (\ref{eq:17}), we note that
\begin{equation}\label{eq:32}
\tilde{C}_{n}^{(k)}(x+y)=\sum_{j=0}^{n}\binom{n}{j}\tilde{C}_{j}^{(k)}(x)p_{n-j}(y),
\end{equation}
where $p_{n}(y)=\frac{1}{Li{f_{k}}(-t)}\tilde{C}_{n}^{(k)}(y)\sim\left(1,e^{t}-1\right)$.\\
By (\ref{eq:23}) and (\ref{eq:24}), we get
\begin{equation}\label{eq:33}
(y)_{n}=p_{n}(y)\sim\left(1,e^{t}-1\right).
\end{equation}
Thus, from (\ref{eq:32}) and (\ref{eq:33}), we have
\begin{equation}\label{eq:34}
\tilde{C}_{n}^{(k)}(x+y)=\sum_{j=0}^{n}\binom{n}{j}\tilde{C}_{j}^{(k)}(x)(y)_{n-j}.
\end{equation}
By (\ref{eq:13}), (\ref{eq:15}) and (\ref{eq:20}), we get
\begin{equation*}
\tilde{C}_{n}^{(k)}(x+1)-\tilde{C}_{n}^{(k)}(x)=\left(e^{t}-1\right)\tilde{C}_{n}^{(k)}(x)=n\tilde{C}_{n-1}^{(k)}(x).
\end{equation*}
For $s_{n}(x)\sim\left(g(t),f(t)\right)$, the recurrence formula for $s_{n}(x)$ is given by
\begin{equation}\label{eq:35}
s_{n+1}(x)=\left(x-\frac{g'(t)}{g(t)}\right)\frac{1}{f'(t)}s_{n}(x),\,\,\,\,(\text{see}\,\,\lbrack 12\rbrack).
\end{equation}
By (\ref{eq:20}) and (\ref{eq:35}), we get
\begin{align}\label{eq:36}
\tilde{C}_{n+1}^{(k)}(x)&=\left(x-\frac{Li{f_{k}'}(-t)}{Li{f_{k}}(-t)}\right)e^{-t}\tilde{C}_{n}^{(k)}(x)\\
&=x\tilde{C}_{n}^{(k)}(x-1)-e^{-t}\frac{Li{f_{k}'}(-t)}{Li{f_{k}}(-t)}\tilde{C}_{n}^{(k)}(x).\nonumber
\end{align}
We observe that
\begin{align}\label{eq:37}
\frac{Li{f_{k}'}(-t)}{Li{f_{k}}(-t)}\tilde{C}_{n}^{(k)}(x)&=Li{f_{k}'}(-t)\frac{1}{Li{f_{k}}(-t)}\tilde{C}_{n}^{(k)}(x)=Li{f_{k}'}(-t)(x)_{n}\\
&=\sum_{l=0}^{n}S_{1}(n,l)Li{f_{k}'}(-t)x^{l}\nonumber\\
&=\sum_{l=0}^{n}S_{1}(n,l)\sum_{m=0}^{l}\frac{(-1)^{m}\binom{l}{m}}{(m+2)^{k}}x^{l-m}\nonumber\\
&=\sum_{j=0}^{n}\left\{\sum_{l=j}^{n}\frac{(-1)^{l-j}\binom{l}{j}}{(l-j+2)^{k}}S_{1}(n,l)\right\}x^{j}.\nonumber
\end{align}
Therefore, by (\ref{eq:36}) and (\ref{eq:37}), we obtain the following theorem.

\begin{thm}\label{eq:thm2}
For $n\geq 0$, we have
\begin{equation*}
\tilde{C}_{n+1}^{(k)}(x)=x\tilde{C}_{n}^{(k)}(x-1)-\sum_{j=0}^{n}\left\{\sum_{l=j}^{n}S_{1}(n,l)\frac{(-1)^{l-j}}{(l-j+2)^{k}}\binom{l}{j}\right\}(x-1)^{j}.
\end{equation*}
\end{thm}

\noindent From (\ref{eq:10}), we note that
\begin{align}\label{eq:38}
\tilde{C}_{n}^{(k)}(y)&=\left\langle\sum_{l=0}^{\infty}\tilde{C}_{l}^{(k)}(y)\frac{t^{l}}{l!}\Bigg\vert x^{n}\right\rangle=\left\langle Li{f_{k}}\left(-\log{(1+t)}\right)(1+t)^{y}\vert x^{n}\right\rangle\\
&=\left\langle Li{f_{k}}\left(-\log{(1+t)}\right)(1+t)^{y}\vert xx^{n-1}\right\rangle\nonumber\\
&=\left\langle \partial_{t}\left(Li{f_{k}}\left(-\log{(1+t)}\right)(1+t)^{y}\right)\vert x^{n-1}\right\rangle\nonumber\\
&=\left\langle \partial_{t}\left(Li{f_{k}}\left(-\log{(1+t)}\right)\right)(1+t)^{y}\vert x^{n-1}\right\rangle\nonumber\\
&\quad +\left\langle Li{f_{k}}\left(-\log{(1+t)}\right)\partial_{t}(1+t)^{y}\vert x^{n-1}\right\rangle\nonumber\\
&=\left\langle\partial_{t}\left(Li{f_{k}}\left(-\log{(1+t)}\right)\right)(1+t)^{y}\vert x^{n-1}\right\rangle +y\tilde{C}_{n-1}^{(k)}(y-1),\nonumber
\end{align}
where $\partial_{t}f(t)=\frac{df(t)}{dt}$.\\
It is easy to show that
\begin{equation}\label{eq:39}
tLi{f_{k}'}(t)=Li{f_{k-1}}(t)-Li{f_{k}}(t).
\end{equation}
Thus, by (\ref{eq:39}), we get
\begin{equation}\label{eq:40}
Li{f_{k}'}(t)=\frac{Li{f_{k-1}}(t)-Li{f_{k}}(t)}{t}.
\end{equation}
By (\ref{eq:38}) and (\ref{eq:40}), we see that
\begin{align}\label{eq:41}
\tilde{C}_{n}^{(k)}(y)&=y\tilde{C}_{n-1}^{(k)}(y-1)\\
&\quad+\left\langle\frac{Li{f_{k-1}}\left(-\log{(1+t)}\right)-Li{f_{k}}\left(-\log{(1+t)}\right)}{(1+t)\log{(1+t)}}(1+t)^{y}\Bigg\vert x^{n-1}\right\rangle\nonumber\\
&=y\tilde{C}_{n-1}^{(k)}(y-1)\nonumber\\
&\quad+\left\langle\frac{Li{f_{k-1}}\left(-\log{(1+t)}\right)-Li{f_{k}}\left(-\log{(1+t)}\right)}{t(1+t)}(1+t)^{y}\Bigg\vert \frac{t}{\log{(1+t)}}x^{n-1}\right\rangle.\nonumber
\end{align}
From (\ref{eq:1}), (\ref{eq:5}) and (\ref{eq:40}), we note that
\begin{align}\label{eq:42}
\tilde{C}_{n}^{(k)}(y)&=y\tilde{C}_{n-1}^{(k)}(y-1)+\sum_{l=0}^{n-1}\frac{B_{l}^{(l)}(1)}{l!}(n-1)_{l}\\
&\quad\times\left\langle\frac{Li{f_{k-1}}\left(-\log{(1+t)}\right)-Li{f_{k}}\left(-\log{(1+t)}\right)}{t}(1+t)^{y-1}\Bigg\vert x^{n-l-1}\right\rangle\nonumber
\end{align}
\begin{align*}
&=y\tilde{C}_{n-1}^{(k)}(y-1)+\sum_{l=0}^{n-1}\frac{B_{l}^{(l)}(1)}{l!}(n-1)_{l}\\
&\quad\times\left\langle\frac{Li{f_{k-1}}\left(-\log{(1+t)}\right)-Li{f_{k}}\left(-\log{(1+t)}\right)}{t}(1+t)^{y-1}\Bigg\vert t\frac{x^{n-l}}{n-l}\right\rangle\\
&=y\tilde{C}_{n-1}^{(k)}(y-1)+\sum_{l=0}^{n-1}\binom{n-1}{l}\frac{B_{l}^{(l)}(1)}{n-l}\left\{\tilde{C}_{n-l}^{(k-1)}(y-1)-\tilde{C}_{n-l}^{(k)}(y-1)\right\}\\
&=y\tilde{C}_{n-1}^{(k)}(y-1)+\frac{1}{n}\sum_{l=0}^{n-1}\binom{n}{l}B_{l}^{(l)}(1)\left\{\tilde{C}_{n-l}^{(k-1)}(y-1)-\tilde{C}_{n-l}^{(k)}(y-1)\right\}.
\end{align*}
It is not difficult to show that $\tilde{C}_{0}^{(k)}(y-1)=\tilde{C}_{0}^{(k-1)}(y-1)$. Therefore, by (\ref{eq:42}), we obtain the following theorem.

\begin{thm}\label{eq:thm3}
For $n\geq 1$, we have
\begin{align*}
\tilde{C}_{n}^{(k)}(x)=x\tilde{C}_{n-1}^{(k)}(x-1)+\frac{1}{n}\sum_{l=0}^{n}\binom{n}{l}B_{l}^{(l)}(1)\left\{\tilde{C}_{n-l}^{(k-1)}(x-1)-\tilde{C}_{n-l}^{(k)}(x-1)\right\}.
\end{align*}
\end{thm}

\noindent For $n\geq m\geq 1$, we compute
\begin{equation*}
\left\langle\left(\log{(1+t)}\right)^{m}Li{f_{k}}\left(-\log{(1+t)}\right)\vert x^{n}\right\rangle
\end{equation*}
in two different ways.\\
On the one hand,
\begin{align}\label{eq:43}
&\left\langle\left(\log{(1+t)}\right)^{m}Li{f_{k}}\left(-\log{(1+t)}\right)\vert x^{n}\right\rangle\\
&=\left\langle Li{f_{k}}\left(-\log{(1+t)}\right)\Bigg\vert\sum_{l=0}^{\infty}\frac{m!}{(l+m)!}S_{1}(l+m,m)t^{l+m}x^{n}\right\rangle\nonumber\\
&=\sum_{l=0}^{n-m}\frac{m!}{(l+m)!}S_{1}(l+m,m)(n)_{l+m}\left\langle Li{f_{k}}\left(-\log{(1+t)}\right)\big\vert x^{n-l-m}\right\rangle\nonumber\\
&=\sum_{l=0}^{n-m}m!\binom{n}{l+m}S_{1}(l+m,m)\tilde{C}_{n-l-m}^{(k)}.\nonumber
\end{align}
On the other hand, we get
\begin{align}\label{eq:44}
&\left\langle\left(\log{(1+t)}\right)^{m}Li{f_{k}}\left(-\log{(1+t)}\right)\vert x^{n}\right\rangle\\
&=\left\langle\left(\log{(1+t)}\right)^{m}Li{f_{k}}\left(-\log{(1+t)}\right)\vert xx^{n-1}\right\rangle\nonumber\\
&=\left\langle\partial_{t}\left(\left(\log{(1+t)}\right)^{m}Li{f_{k}}\left(-\log{(1+t)}\right)\right) \vert x^{n-1}\right\rangle.\nonumber
\end{align}
Now, we observe that
\begin{align}\label{eq:45}
&\partial_{t}\left(\left(\log{(1+t)}\right)^{m}Li{f_{k}}\left(-\log{(1+t)}\right)\right)\\
&=m\left(\log{(1+t)}\right)^{m-1}\frac{1}{1+t}Li{f_{k}}\left(-\log{(1+t)}\right)\nonumber\\
&\quad+\left(\log{(1+t)}\right)^{m}\frac{Li{f_{k-1}} \left(-\log{(1+t)}\right)-Li{f_{k}}\left(-\log{(1+t)}\right)}{(1+t)\log{(1+t)}}\nonumber\\
&=\left(\log{(1+t)}\right)^{m-1}(1+t)^{-1}\nonumber\\
&\quad\times\left\{mLi{f_{k}}\left(-\log{(1+t)}\right)+Li{f_{k-1}}\left(-\log{(1+t)}\right)-Li{f_{k}}\left(-\log{(1+t)}\right)\right\}.\nonumber
\end{align}
By (\ref{eq:44}) and (\ref{eq:45}), we get
\begin{align}\label{eq:46}
&\left\langle\left(\log{(1+t)}\right)^{m}Li{f_{k}}\left(-\log{(1+t)}\right)\big\vert x^{n}\right\rangle\\
&=\sum_{l=0}^{n-m}\frac{(m-1)!}{(l+m-1)!}S_{1}(l+m-1,m-1)\nonumber\\
&\quad\times\{ (m-1)\left\langle Li{f_{k}}\left(-\log{(1+t)}\right)(1+t)^{-1}\big\vert t^{l+m-1}x^{n-1}\right\rangle\nonumber\\
&\quad +\left\langle Li{f_{k-1}}\left(-\log{(1+t)}\right)(1+t)^{-1}\vert t^{l+m-1}x^{n-1}\right\rangle\}\nonumber\\
&=(m-1)\sum_{l=0}^{n-m}\frac{(m-1)!}{(l+m-1)!}S_{1}(l+m-1,m-1)(n-1)_{l+m-1}\nonumber\\
&\quad\times\left\langle Li{f_{k}}\left(-\log{(1+t)}\right)(1+t)^{-1}\vert x^{n-m-l}\right\rangle\nonumber\\
&\quad+\sum_{l=0}^{n-m}\frac{(m-1)!}{(l+m-1)!}S_{1}(l+m-1,m-1)(n-1)_{l+m-1}\nonumber\\
&\quad\times\left\langle Li{f_{k-1}}\left(-\log{(1+t)}\right)(1+t)^{-1}\vert x^{n-m-l}\right\rangle\nonumber\\
&=\sum_{l=0}^{n-m}(m-1)!\binom{n-1}{l+m-1}S_{1}(l+m-1,m-1)\nonumber\\
&\quad\times\left\{(m-1)\tilde{C}_{n-l-m}^{(k)}(-1)+\tilde{C}_{n-l-m}^{(k-1)}(-1)\right\}.\nonumber
\end{align}
Therefore, by (\ref{eq:43}) and (\ref{eq:46}), we obtain the following theorem.

\begin{thm}\label{eq:thm4}
For $n\geq m\geq 1$, we have
\begin{align*}
&\sum_{l=0}^{n-m}m!\binom{n}{l+m}S_{1}(l+m,m)\tilde{C}_{n-l-m}^{(k)}\\
&=\sum_{l=0}^{n-m}(m-1)!\binom{n-1}{l+m-1}S_{1}(l+m-1,m-1)\\
&\quad\times\left\{(m-1)\tilde{C}_{n-l-m}^{(k)}(-1)+\tilde{C}_{n-l-m}^{(k-1)}(-1)\right\}.
\end{align*}
In particular, if we take $m=1$, then we get
\begin{equation*}
\tilde{C}_{n}^{(k-1)}(-1)=\sum_{l=0}^{n-1}(-1)^{l}l!\binom{n}{l+1}\tilde{C}_{n-l-1}^{(k)}.
\end{equation*}
\end{thm}

\noindent$\large{\mathbf{Remark}}$. For $s_{n}(x)\sim\left(g(t),f(t)\right)$, it is known that
\begin{equation}\label{eq:47}
\frac{d}{dx}s_{n}(x)=\sum_{l=0}^{n-1}\binom{n}{l}\left\langle\bar{f}(t)\big\vert x^{n-l}\right\rangle s_{l}(x).
\end{equation}
By (\ref{eq:20}) and (\ref{eq:47}), we easily show that
\begin{equation*}
\frac{d}{dx}\tilde{C}_{n}^{(k)}(x)=(-1)^{n}n!\sum_{l=0}^{n-1}\frac{(-1)^{l-1}}{(n-l)l!}\tilde{C}_{l}^{(k)}(x).
\end{equation*}
Let us consider the following two Sheffer sequences:
\begin{equation}\label{eq:48}
\tilde{C}_{n}^{(k)}(x)\sim\left(\frac{1}{Li{f_{k}}(-t)},e^{t}-1\right),
\end{equation}
and
\begin{equation*}
B_{n}^{(r)}(x)\sim\left(\left(\frac{e^{t}-1}{t}\right)^{r},t\right).
\end{equation*}
Suppose that
\begin{equation}\label{eq:49}
\tilde{C}_{n}^{(k)}(x)=\sum_{m=0}^{n}C_{n,m}B_{m}^{(r)}(x).
\end{equation}
By (\ref{eq:19}), we see that
\begin{align}\label{eq:50}
C_{n,m}&=\frac{1}{m!}\left\langle\frac{\left(\frac{t}{\log{(1+t)}}\right)^{r}}{\frac{1}{Li{f_{k}}\left(-\log{(1+t)}\right)}}\left(\log{(1+t)}\right)^{m}\Bigg\vert x^{n}\right\rangle\\
&=\frac{1}{m!}\left\langle Li{f_{k}}\left(-\log{(1+t)}\right)\left(\frac{t}{\log{(1+t)}}\right)^{r}\left(\log{(1+t)}\right)^{m}\Bigg\vert x^{n}\right\rangle\nonumber\\
&=\frac{1}{m!}\sum_{l=0}^{n-m}\frac{m!}{(l+m)!}S_{1}(l+m,m)(n)_{l+m}\nonumber\\
&\quad\times\left\langle Li{f_{k}}\left(-\log{(1+t)}\right)\left(\frac{t}{\log{(1+t)}}\right)^{r}\Bigg\vert x^{n-l-m}\right\rangle\nonumber\\
&=\frac{1}{m!}\sum_{l=0}^{n-m}\frac{m!}{(l+m)!}S_{1}(l+m,m)(n)_{l+m}\sum_{a=0}^{n-l-m}B_{a}^{(a-r+1)}\frac{1}{a!}\nonumber\\
&\quad\times\left\langle Li{f_{k}}\left(-\log{(1+t)}\right)\big\vert t^{a}x^{n-l-m}\right\rangle\nonumber\\
&=\sum_{l=0}^{n-m}\binom{n}{l+m}S_{1}(l+m,m)\sum_{a=0}^{n-l-m}B_{a}^{(a-r+1)}\frac{(n-l-m)_{a}}{a!}\nonumber\\
&\quad\times\left\langle Li{f_{k}}\left(-\log{(1+t)}\right)\big\vert x^{n-l-m-a}\right\rangle\nonumber\\
&=\sum_{l=0}^{n-m}\sum_{a=0}^{n-l-m}\binom{n}{l+m}\binom{n-m-l}{a}S_{1}(l+m,m)B_{a}^{(a-r+1)}(1)\tilde{C}_{n-l-m-a}^{(k)}.\nonumber
\end{align}
Therefore, by (\ref{eq:49}) and (\ref{eq:50}), we obtain the following theorem.

\begin{thm}\label{eq:thm5}
For $n\geq 0$, we have
\begin{align*}
\tilde{C}_{n}^{(k)}(x)&=\sum_{m=0}^{n}\Bigg\{\sum_{l=0}^{n-m}\sum_{a=0}^{n-m-l}\binom{n}{l+m}\binom{n-m-l}{a}S_{1}(l+m,m)B_{a}^{(a-r+1)}(1)\\
&\quad\times\tilde{C}_{n-m-l-a}^{(k)}\Bigg\}B_{m}^{(r)}(x).
\end{align*}
\end{thm}

\noindent$\large{\mathbf{Remark}}$. The Narumi polynomials of order $a$ are defined by the generating function to be
\begin{equation}\label{eq:51}
\sum_{k=0}^{\infty}\frac{N_{k}^{(a)}(x)}{k!}t^{k}=\left(\frac{t}{\log{(1+t)}}\right)^{-a}(1+t)^{x}.
\end{equation}
Indeed, $N_{a}^{(k)}(x)=B_{k}^{(k+a+1)}(x+1)$, $N_{k}^{(a)}(x)\sim\left(\left(\frac{e^{t}-1}{t}\right)^{a},e^{t}-1\right)$.\\
By (\ref{eq:50}) and (\ref{eq:51}), we get
\begin{equation}\label{eq:52}
C_{n,m}=\sum_{l=0}^{n-m}\sum_{a=0}^{n-m-l}\binom{n}{l+m}\binom{n-l-m}{a}S_{1}(l+m,m)N_{a}^{(-r)}\tilde{C}_{n-l-m-a}^{(k)}.
\end{equation}
From (\ref{eq:49}) and (\ref{eq:52}), we have
\begin{align}\label{eq:53}
&\tilde{C}_{n}^{(k)}(x)\\
&=\sum_{m=0}^{n}\left\{\sum_{l=0}^{n-m}\sum_{a=0}^{n-m-l}\binom{n}{l+m}\binom{n-l-m}{a}S_{1}(l+m,m)N_{a}^{(-r)}\tilde{C}_{n-l-m-a}^{(k)}\right\}B_{m}^{(r)}(x).\nonumber
\end{align}
By (\ref{eq:1}), we easily show that
\begin{align}\label{eq:54}
C_{n,m}&=\sum_{l=0}^{n-m}\sum_{a=0}^{n-m-l}\sum_{a_{1}+\cdots+a_{r}=a}^{}\binom{n}{l+m}\binom{n-l-m}{a}\binom{a}{a_{1},\cdots,a_{r}}\\
&\quad\times S_{1}(l+m,m)b_{a_{1}}\cdots b_{a_{r}}\tilde{C}_{n-m-l-a}^{(k)}.\nonumber
\end{align}
From (\ref{eq:49}) and (\ref{eq:54}), we can derive the following equation:
\begin{align}\label{eq:55}
\tilde{C}_{n}^{(k)}(x)&=\sum_{m=0}^{n}\Bigg\{\sum_{l=0}^{n-m}\sum_{a=0}^{n-m-l}\sum_{a_{1}+\cdots +a_{r}=a}^{}\binom{n}{l+m}\binom{n-l-m}{a}\binom{a}{a_{1},\cdots,a_{r}}\\
&\quad\times S_{1}(l+m,m)\left(\prod_{i=1}^{r}b_{a_{i}}\right)\tilde{C}_{n-m-l-a}^{(k)}\Bigg\}B_{m}^{(r)}(x).\nonumber
\end{align}
For $\tilde{C}_{n}^{(k)}(x)\sim\left(\frac{1}{Li{f_{k}}(-t)}, e^{t}-1\right)$, $H_{n}^{(r)}(x\vert\lambda)\sim\left(\left(\frac{e^{t}-\lambda}{1-\lambda}\right)^{r},t\right)$, $(r\geq 0)$, let
\begin{equation}\label{eq:56}
\tilde{C}_{n}^{(k)}(x)=\sum_{m=0}^{n}C_{n,m}H_{m}^{(r)}(x\vert\lambda),
\end{equation}
where, by (\ref{eq:19}), we get
\begin{align}\label{eq:57}
C_{n,m}&=\frac{1}{m!(1-\lambda)^{r}}\left\langle Li{f_{k}}\left(-\log{(1+t)}\right)(1+t-\lambda)^{r}\big\vert\left(\log{(1+t)}\right)^{m}x^{n}\right\rangle\\
&=\frac{1}{m!(1-\lambda)^{r}}\sum_{l=0}^{n-m}\frac{m!}{(l+m)!}S_{1}(l+m,m)(n)_{l+m}\nonumber\\
&\quad\times\left\langle Li{f_{k}}\left(-\log{(1+t)}\right)(1+t-\lambda)^{r}\big\vert x^{n-l-m}\right\rangle.\nonumber
\end{align}
We observe that
\begin{align}\label{eq:58}
&\left\langle Li{f_{k}}\left(-\log{(1+t)}\right)(1+t-\lambda)^{r}\big\vert x^{n-l-m}\right\rangle\\
&=\sum_{a=0}^{r}\binom{r}{a}(1-\lambda)^{r-a}\left\langle  Li{f_{k}}\left(-\log{(1+t)}\right)\big\vert t^{a}x^{n-l-m}\right\rangle\nonumber\\
&=\sum_{a=0}^{r}\binom{r}{a}(1-\lambda)^{r-a}(n-m-l)_{a}\left\langle Li{f_{k}}\left(-\log{(1+t)}\right)\big\vert x^{n-l-m-a}\right\rangle\nonumber\\
&=\sum_{a=0}^{r}\binom{r}{a}(1-\lambda)^{r-a}(n-m-l)_{a}\tilde{C}_{n-l-m-a}^{(k)}.\nonumber
\end{align}
Thus, by (\ref{eq:57}) and (\ref{eq:58}), we get
\begin{equation}\label{eq:59}
C_{n,m}=\sum_{l=0}^{n-m}\sum_{a=0}^{r}\binom{n}{l+m}\binom{r}{a}(n-m-l)_{a}(1-\lambda)^{-a}S_{1}(l+m,m)\tilde{C}_{n-m-l-a}^{(k)}.
\end{equation}
Therefore, by (\ref{eq:56}) and (\ref{eq:59}), we obtain the following theorem.

\begin{thm}\label{eq:thm6}
For $n\geq 0$, we have
\begin{align*}
\tilde{C}_{n}^{(k)}(x)&=\sum_{m=0}^{n}\Bigg\{\sum_{l=0}^{n-m}\sum_{a=0}^{r}\binom{n}{l+m}\binom{r}{a}(n-m-l)_{a}(1-\lambda)^{-a}S_{1}(l+m,m)\\
&\quad\times\tilde{C}_{n-m-l-a}^{(k)}\Bigg\}H_{m}^{(r)}(x\vert\lambda).
\end{align*}
\end{thm}

\noindent For $\tilde{C}_{n}^{(k)}(x)\sim\left(\frac{1}{Li{f_{k}}(-t)},e^{t}-1\right)$, and $(x)_{n}\sim\left(1,e^{t}-1\right)$, let us assume that
\begin{equation}\label{eq:60}
\tilde{C}_{n}^{(k)}(x)=\sum_{m=0}^{n}C_{n,m}(x)_{m}.
\end{equation}
From (\ref{eq:19}), we note that
\begin{align}\label{eq:61}
C_{n,m}&=\frac{1}{m!}\left\langle Li{f_{k}}\left(-\log{(1+t)}\right)t^{m}\big\vert x^{n}\right\rangle\\
&=\frac{1}{m!}\left\langle Li{f_{k}}\left(-\log{(1+t)}\right)\big\vert t^{m}x^{n}\right\rangle\nonumber\\
&=\binom{n}{m}\left\langle Li{f_{k}}\left(-\log{(1+t)}\right)\big\vert x^{n-m}\right\rangle\nonumber\\
&=\binom{n}{m}\tilde{C}_{n-m}^{(k)}\nonumber.
\end{align}
Therefore, by (\ref{eq:60}) and (\ref{eq:61}), we obtain the following theorem.

\begin{thm}\label{eq:thm7}
For $n\geq 0$, we have
\begin{equation*}
\tilde{C}_{n}^{(k)}(x)=\sum_{m=0}^{n}\binom{n}{m}\tilde{C}_{n-m}^{(k)}(x)_{m}.
\end{equation*}
\end{thm}


\bigskip
ACKNOWLEDGEMENTS. This work was supported by the National Research Foundation of Korea(NRF) grant funded by the Korea government(MOE)\\
(No.2012R1A1A2003786 ).
\bigskip

\noindent
\author{Department of Mathematics, Sogang University, Seoul 121-742, Republic of Korea
\\e-mail: dskim@sogang.ac.kr}\\
\\
\author{Department of Mathematics, Kwangwoon University, Seoul 139-701, Republic of Korea
\\e-mail: tkkim@kw.ac.kr}
\end{document}